# An Environmentally Sustainable Closed-Loop Supply Chain Network Design under Uncertainty: Application of Optimization

Md. Mohsin Ahmed[1], S. M. Salauddin Iqbal[2], Tazrin Jahan Priyanka[3], Mohammad Arani[4][0000-0002-1712-067X], Mohsen Momenitabar[5][0000-0003-2568-1781], and Md Mashum Billal[6][0000-0001-5557-3613]

[1] Department of Industrial & Production Engineering, Rajshahi University of Engineering & Technology (RUET), Rajshahi-6204, Bangladesh
durjoymohsin@gmail.com

[2] Department of Industrial & Production Engineering, Rajshahi University of Engineering & Technology (RUET), Rajshahi-6204, Bangladesh
salauddin11488@gmail.com

[3] Department of Industrial & Production Engineering, Rajshahi University of Engineering & Technology (RUET), Rajshahi-6204, Bangladesh
priyanka.ruet12@gmail.com

[4] Department of Systems Engineering, University of Arkansas at Little Rock, USA
mxarani@ualr.edu

[5] Department of Transportation, Logistics, and Finance, North Dakota State University, Fargo, USA
mohsen.momenitabar@ndsu.edu

[6] Department of Mechanical Engineering (Engineering Management), University of Alberta, Edmonton, Canada
mdmashum@ualberta.ca

**Abstract.** Newly, the rates of energy and material consumption to augment industrial production are substantially high, thus the environmentally sustainable industrial development has emerged as the main issue of either developed or developing countries. A novel approach to supply chain management is proposed to maintain economic growth along with environmentally friendly concerns for the design of the supply chain network. In this paper, a new green supply chain design approach has been suggested to maintain the financial virtue accompanying the environmental factors that required to be mitigated the negative effect of rapid industrial development on the environment. This approach has been suggested a multi-objective mathematical model minimizing the total costs and CO2 emissions for establishing an environmentally sustainable closed-loop supply chain. Two optimization methods are used namely Epsilon Constraint Method, and Genetic Algorithm Optimization Method. The results of the two mentioned methods have been compared and illustrated their effectiveness. The outcome of the analysis is approved to verify the accuracy of the proposed model to deal with financial and environmental issues concurrently.

# 1    Introduction

The sustainable supply chain has been one of the most uprising topics in the current decade as the environment is polluted at an alarming rate due to the functioning of organizations and factories. Sustainable supply chain (SSC) is, thus, a management system that integrates environmental concerns with supply chain functioning such as product design, raw material sourcing, manufacturing, assembly, transportation of final goods, and controlling the disposal and recycling of products when their life has ended (Aldemir et al. 2018). The use of the green supply chain (GSC) is generally originated from two common interrelated factors the amount of carbon emission and the cost of claiming a sustainable supply chain. Environmental issues have already absorbed the attention of the developed countries; federal law and customer pressure have increased to follow SSC and considering environmentally friendly production, consequently (Niranjan et al. 2019).

For green uplift of the total supply chain network, carbon emission at every entity must be taken into consideration. Focused on the green effect on a closed-looped supply chain and attempted to abate $CO_2$ emissions, furthermore, motivate the customers to utilize the recyclable product. In our research, we have formulated a mathematical model to optimize both the carbon emission and the total cost. Here, multi-objective optimization has been identified using two mathematical approaches both the exact solution method and meta-heuristic approach. A closed-loop supply chain network has also been incorporated with reverse flow. This closed-loop mainly deals with the end life of products and unused raw materials that otherwise it can create more wastage and degradation of the environment. We have also considered mixed-integer linear programming for more comprehensive supply chain network evaluation and thus considering green development. (Konyalıoğlu and Zafeirakopoulos 2020; Yang et al. 2018)

In our study, the only natural environment is being considered. The natural environment is the interaction between biology, ecology, and other forces. To maintain a safe environment this interaction must be maintained and balanced. Though industrialization that increases the way of living standard, it also has a huge escapade effect on the environment; manufacturing emits a lot of wastage and greenhouse gases that affect environmental balance. However, this topic remained unnoticed up to the 20th century. Concern for the pollution developed a concern for the environment thus leading us to save ourselves and the earth. People are now more concerned than ever before as anthropogenic activity already damage the balance of ecology and the environment. As a result, several agreements have been signed on compliance with environmental objectives. Some most important protocol and conventions are:

- Montreal protocol: It was signed to save the ozone layer and signed in 1987. The Montreal Protocol stipulated that production and use of ozone-depleting chemicals would be phased out by 2000. The primary purpose was to reduce CFCs.

- Kyoto protocol: This was a consensus reached at the 1992 Earth Summit in Rio de Janeiro under the United Nations Framework Convention on Climate Change (UNFCCC). The primary goal was to maintain a small enough degree of reduction of greenhouse gas emissions in the environment to avoid harmful anthropogenic involvement with the climate system. The Protocol has been ratified by 154 countries.

## 2    Literature Review

The sustainable supply chain is a set of procedures performed to monitor and facilitate the environmental and economic performance indices by distributing required resources, material, goods, and allocating organizational responsibilities and tasks (Kim et al. 2011). GSC may also be used as reducing and preferably alleviating the adverse impact of the supply chain network on the natural ecosystem (Andiç et al. 2012). In our research on the green supply chain, the literature concerned about finding the sectors and subsectors where minimization of carbon could be justified and regulate the emission of carbon and the cost of the total supply chain. In recent years mostly in the 21st century, it is the most stated topics and concern has given on environment-friendly production and supply chain. So, there are several innovative studies have been done on this concern in recent times.

Green supply chain (GSC) refers to a strategic approach considering environmental indices to design the supply chain processes. This includes a group of sub-processes like product design, raw material sourcing, manufacturing, transportation, and lifecycle management. According to Wee et al. (2011), GSC is incorporating environmental problems into supply chain management, including packaging design, content procurement, and quality, production operations, distribution of the finished product to customers, and end-of-life management of greening products. GSC thus just not merely a concept of green manufacturing, instead, it is the integration of all operations of the supply chain that needs to be environment friendly. In other words, GSC can be narrowly defined as a significant new model for businesses to achieve income and market share goals by reducing their environmental risks and impacts while increasing their environmental performance. GSC may be viewed as responsible for each segment of the overall supply chain network. Some sectors are as follows:

- *Assembly:* Carbon is emitted while assembling and manufacturing a product.
- *Disassembly:* While performing disassembly after products useful life, the operation requires fossil fuels or electricity.
- *Transportation:* This produces a huge amount of carbon emission. Transportation means like truck, launch, ship, forklift, etc. use a lot of fossil fuels like diesel, petrol, coil, etc. and burning of this produces carbon and other greenhouse gases.
- *Handling:* Handling of raw materials, machinery, equipment, etc. require some electro-mechanical devices and power.
- *Remanufacturing:* In a closed-loop supply chain, after the useful life of a product is over, it needs to be remanufactured thus reducing its impact on the environment.

Che et al. (2017) examined a multi-floor facility layout problem considering bi-objective formulation. To illustrate this problem, they modeled a bi-objective mixed-integer non-linear programming model. The main objective was to minimize the total material handling costs and the total engaged room area. They have used the ε-constraint approach to evaluate the solution and thus to figure out the optimum location. Sinha and Chaturvedi (2018), investigated two important factors in production planning for green development. They studied two objectives of carbon emission and energy consumption and offered manufacturing planning through available and new process routes. They suggested a graphical method according to the perceptions of pinch analysis. They also showed that the proposed model will reduce the extra burden of $CO_2$ emission and energy consumption while meeting the demand.

Talaei et al. (2016), worked on reducing carbon emission in the overall supply chain network for an environmental issue. They also considered the total cost of the supply chain. They used robust fuzzy programming to investigate this issue. They considered the ε-constraint approach and numerical model to illustrate and solve this optimization problem. Nurjanni et al. (2017), investigated the traditional supply chain management and showed that with increasing levels of industrialization and globalization, the supply chain should be more robust and efficient. They focused on green development and proposed a new supply chain network that considered environmental and financial issues. They were proposing a concept of multi-objective optimization to reduce costs and carbon emissions. Three scalarization methods were used to test the optimization process, namely, the weighted sum method, and weighted Tchebycheff, and increased weighted Tchebycheff to illustrate the trade-offs between costs and environmental concerns. Cheng et al. (2018), worked on the transit system and aimed at public transit systems with optimum designs and operating plans that can mitigate both total costs and greenhouse gas (GHG) emissions. This study investigated the potential emission impact when optimizing the total costs of the transit system. They used a continuum approximation (CA) method to show the relationship of the network. They analyzed the optimal total costs and emissions and the interaction between them. This research showed a significant reduction in GHG emissions that can be achieved when total costs are reduced simultaneously. Additionally another instance of GHG emission control could be find Dehdari Ebrahimi and Momeni Tabar (2017). Lahnaoui et al. (2018) proposed a method to identify the minimum cost of constructing a hydrogen infrastructure using a single objective linear optimization. They only consider the minimization of cost. This work has focused on reducing both hydrogen transportation capital and operational costs. They consider costs associated with setting up inventory and compression facilities as well as the transport network. Hence, the main aim of the study was to create a cost-effective transport network using compressed gas trucks for mobility. They also showed potential hydrogen demand figures in the country, which will rise according to their model from 2030 to 2050. The research also concentrated on whether and at what flow demand for hydrogen would be delivered from the manufacturing nodes into the global distribution network. Nie et al. (2018) worked on optimizing electric power systems with cost minimization and environmental impact reduction. They used a multistage inexact-factorial fuzzy-probability programming method to illustrate their mathematical form. They also investigated quantitative analysis for individual factors and their interaction. Their policy was to analyze individual, uncertainty recognition, and interaction consideration. They showed that the uncertainties influence Qingdao's EPS long-term planning which can be organized to a cleaner and safer model exhibiting renewable energy. They showed the relation between electricity demand and associated cost.

**Table 1.** Table captions should be placed above the tables.

| Author | Net-work | Model | Objec-tive | Solution | Sectors Under Consideration $CO_2$ Emission | | | | |
|--------|----------|-------|------------|----------|------|------|------|------|------|
| | | | | | PRO | HL | TRS | DA | RM |
| Diabat and Simchi-Levi (2009) | F | MILP | MO | Weighted Sum | C | - | C | C | C |
| Franca et al. (2010) | OL | MILP | MO | 6 Sigma | C | - | - | - | - |
| Wang et al. (2011) | OL | MILP | MO | Facility Layout | C | C | C | - | - |
| Abdallah et al. (2012) | OL | MILP | SO | Net. De-sign | C | - | C | - | - |
| Elhedhli and Merrick (2012) | OL | MILP | MO | Lagrange relaxa-tion | C | C | C | - | - |
| Zhang and Liu (2013) | OL | NLP | SO | Game theory I, II | C | - | - | - | - |
| Talaei et al. (2016) | CL | MILP | MO | ε-con-straint | C | - | - | - | - |
| Nurjanni et al. (2017) | CL | MILP | MO | Weighted Sum | C | C | C | C | C |
| Nie et al. (2018) | OL | MILP | MO | Fuzzy-Probabil-ity | C | C | - | - | - |
| Cheng et al. (2018) | CL | MILP | MO | Contin-uum Ap-proxima-tion | C | C | - | - | - |
| Billal and Hossain (2020) | F | MILP | MO | Stochas-tic Opti-mization | C | C | C | - | - |
| This Research | CL | MILP | MO | ε-con-straint & GA | C | C | C | C | C |

NOTE: F = Forward logistics; RL = Reverse Logistics; OL = Open Loop; CL = Closed-Loop; SC = Supply Chain; LP = Linear Programming; PNLP = Progressive Non-Linear Programming; MILP = Mixed Integer Linear Programming; NLP = Non-Linear Programming; SO = Single Objective; MO = Multi Objective; C = Considered, DA=Disassemble, RM= Remanufacturing, TRS= Transport, HL=Handling.

# 3 Problem Definition

The GSC problem mathematical model consists of four types of echelons, namely plants, warehouses, clients, and disassembly centers which are shown in **Fig. 1**. It is assumed to be one product solely. There are some substitutes for transport (i.e. roads, railways, water, etc.) within each connection link. The structure of the product is not always identical because a specified quantity of product units is allotted by such connections and different arrangements are implemented at each stage (or echelon). Three types of products form (new product, product for disposal, and product for demolition) run through the mathematical model through the supply chain network, in which the product motion description, based on product structure, is as follows:

- A new product, which is new or produced, would be forwarded along the network line (factory to warehouse and warehouse to a customer) to meet customers' demand.
- A product to be disposed of is collected from the customer and sent to a recycling center for disassembly.
- A product to be demolished is distributed from the disassembly center to the reproducing plant.

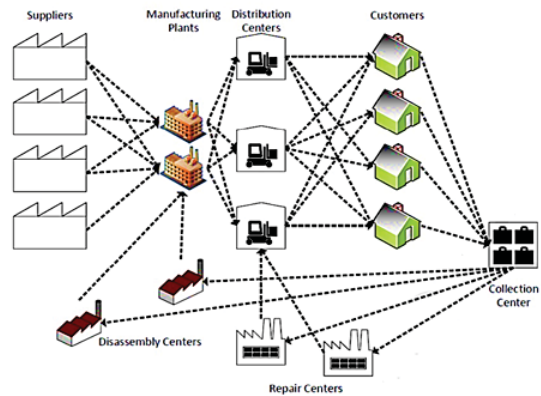

**Fig. 1.** Configuration of the connection between the stages of the supply chain (Alegoz et al. 2020).

### 3.1 Assumptions

1. The demand of all customers is uncertain.
2. Each product unit manufactured at a plant is free of defects.
3. Warehouse reliability is considered.
4. In all warehouses, every client demand is always met by any factory.
5. The total product units to be disposed of that enter a disassembly center is effectively disassembled.
6. There is no restriction regarding the capacity of available transportation options.
7. It is not possible to serve customer demand from intermodal transport options (road, rail, sea, and air).
8. The distances between the network junctions used for the assessment of $CO_2$ emissions are determined by the straight line between them. Transportation rates are used to calculate the characteristics of different types of networks (curviness', impossible link, etc.).

## 3.2 Notations

The following notation is used after observation of the literature and assuming practical situations.

**Table 2.** Parameters.

| Parameters | Definition | Parameters | Definition |
|---|---|---|---|
| $F$ | Potential Factories | $W$ | Potential Warehouses |
| $C$ | Customers | $I$ | Potential disassembly Centers |
| $TF$ | Transportation from factories | $TW$ | Transportation from warehouses |
| $TK$ | transportation from customers | $TI$ | Transportation of disassembly Centers |
| $Q_c$ | Customers' demand | $Ta_{fw}^{tf}$ | Unit Cost of transport from plant to warehouse |
| $Tb_{wc}^{tw}$ | Unit Costs of carriage from a warehouse to a customer | $Tc_{ci}^{tk}$ | Unit Costs of carriage for gathering units of product that needs to be disposed of from customers to disassembly center |
| $Td_{if}^{ti}$ | Unit Costs of carriage from disassembly center to the factory | $L_{fw}^{tf}$ | Factory-to-warehouse transport rate |
| $L_{wc}^{tw}$ | Warehouse-to-customer transportation rate | $L_{ci}^{tk}$ | Transportation rate for the collection of product units to be disposed of from the customer to the disassembly center |
| $L_{if}^{ti}$ | Transportation rate from disassembly center to the factory | $Da_{fw}$ | Distance between factory and warehouse |
| $Db_{wc}$ | Distance between warehouse and customer | $Dc_{ci}$ | Distance between client and disassembly center |
| $Dd_{if}$ | Distance from the disassembly center to the factory | $Ra_f$ | Established plant fixed set-up costs |
| $Rb_w$ | Established warehouse fixed set-up costs | $Rd_i$ | Established disassembly center fixed set-up costs |
| $Ma_f$ | Variable unit costs to produce a product in the factory | $Mb_w$ | Variable unit costs to handle a product in the warehouse |
| $Mc_c$ | Unit variable costs for collecting a disposal product from customers | $Md_i$ | Unit variable cost for disassembling a product to be disposed of in a disassembly center |
| $Mr_f$ | Unit variable cost for remanufacturing a product to be demolished in the factory | $Pa_f$ | Maximum production capacity of the factory |
| $Pb_w$ | Maximum processing capacity of the warehouse | $Pd_i$ | Disassembly Center 's maximum capacity |
| $Pr_f$ | Remanufacturing factory's maximum capacity | $Hd$ | Minimum percentage of the units of the commodity to be disposed of to be obtained from consumers |
| $Hr$ | Minimum percentage of the units of the product to be demolished to be shipped from the disassembly center | $Ga_f$ | Level of $CO_2$ emitted for the development of one unit in a factory |
| $Gb_w$ | Level of $CO_2$ emitted for handling one commodity in a warehouse | $Gd_i$ | Level of $CO_2$ emitted to disassemble one unit of product to be disposed of in a disassembly center |
| $Gr_f$ | Level of $CO_2$ emitted to remanufacture one unit of product to be demolished in a factory | $Gta^{tf}$ | Level of $CO_2$ emitted by shipping a unit of product from a factory to a warehouse for a unit distance |
| $Gtb^{tw}$ | Level of $CO_2$ emitted by shipping a unit of product from a warehouse to a customer for a unit distance | $Gtc^{tk}$ | Level of $CO_2$ emitted by shipping to collect a unit of disposal from customer to disassembly center for a unit distance |
| $Gtd^{ti}$ | Level of $CO_2$ emitted by shipping a unit of product to be demolished from disassembly center to factory for a unit distance | | |

**Table 3.** Decision Variables.

| Variable | Definition | Variable | Definition |
|---|---|---|---|
| $Xa_f$ | 1 if the factory is opened, and 0 in otherwise | $Xb_w$ | 1 if the warehouse is opened, and 0 in otherwise |
| $Xd_i$ | 1 if the disassembly center is opened, and 0 in otherwise | $Ya_{fw}^{tf}$ | Quantity of product units shipped from factory to warehouse, considering the transport options |
| $Yb_{wc}^{tw}$ | Quantity of product units shipped from a warehouse to customer, considering the transport options | $Yc_{ci}^{tk}$ | Quantity of product units to be disposed of shipped from a customer to disassembly center, considering the transport options |
| $Yd_{if}^{tl}$ | Quantity of product units to be demolished shipped from a disassembly center to factory, considering the transport options | $e^{-\lambda t}$ | Reliability factor: where $t$ is time, and $\lambda$ is the failure rate of a warehouse. |

### 3.3 Mathematical Formulation

The first function (2) works for minimizing the total cost. The total cost denotes the total fixed cost, the total variable cost, and the total transportation cost.

$$Min\ f_1 = TC = TFC + TVC + TTC \tag{1}$$

Here,

$$TFC = \sum_{f \in F} Ra_f Xa_f + \sum_{w \in W} Rb_w Xb_w + \sum_{i \in I} Rd_i Xd_i \tag{2}$$

$$TVC = \sum_{f \in F} Ma_f \sum_{w \in W} \sum_{tf \in TF} Ya_{fw}^{tf} + e^{-\lambda t} \sum_{w \in W} Mb_w \sum_{c \in C} \sum_{tw \in TW} Yb_{wc}^{tw} + \sum_{c \in C} Mc_c \sum_{i \in I} \sum_{tk \in TK} Yc_{ci}^{tk} + \sum_{i \in I} Md_i \sum_{c \in C} \sum_{tk \in TK} Yc_{ci}^{tk} + \sum_{f \in F} Mr_f \sum_{i \in I} \sum_{tl \in Tl} Yd_{if}^{tl} \tag{3}$$

$$TTC = \sum_{f \in F} \sum_{w \in W} \sum_{tf \in TF} Ta_{fw}^{tf} Ya_{fw}^{tf} + e^{-\lambda t} \sum_{w \in W} \sum_{c \in C} \sum_{tw \in TW} Tb_{wc}^{tw} Yb_{wc}^{tw} + \sum_{c \in C} \sum_{i \in I} \sum_{tk \in TK} Tc_{ci}^{tk} Yc_{ci}^{tk} + \sum_{i \in I} \sum_{f \in F} \sum_{tl \in Tl} Td_{if}^{tl} Yd_{if}^{tl} \tag{4}$$

The second function (6) works for minimizing the total emission of $CO_2$. The total emission of $CO_2$ denotes the emission of production, assembly, handling, disassembly, remanufacturing, and transportation.

$$Min\ f_2 = TE = EP + EH + ED + ER + ET \tag{5}$$

Here,

$$EP = \sum_{f \in F} Ga_f \sum_{w \in W} \sum_{tf \in TF} Ya_{fw}^{tf} \tag{6}$$

$$EA = \sum_{f \in F} Gc_f \sum_{w \in W} \sum_{tf \in TF} Ya_{fw}^{tf} \tag{7}$$

$$EH = e^{-\lambda t} \sum_{w \in W} Gb_w \sum_{c \in C} \sum_{tw \in TW} Yb_{wc}^{tw} \tag{8}$$

$$ED = \sum_{i \in I} Gd_i \sum_{c \in C} \sum_{tk \in TK} Yc_{ci}^{tk} \tag{9}$$

$$ER = \sum_{f \in F} Gr_f \sum_{i \in I} \sum_{ti \in TI} Yd_{if}^{ti} \tag{10}$$

$$ET = \sum_{tf \in TF} Gta^{tf} \sum_{f \in F} \sum_{w \in W} Ya_{fw}^{tf} Da_{fw} L_{fw}^{tf} + e^{-\lambda t} \sum_{tw \in TW} Gtb^{tw} + \sum_{w \in W} \sum_{c \in C} Yb_{wc}^{tw} Db_{wc} L_{wc}^{tw}$$
$$\sum_{tk \in TK} Gtc^{tk} \sum_{c \in C} \sum_{i \in I} Yc_{ci}^{tk} Dc_{ci} L_{ci}^{tk} + \sum_{ti \in TI} Gtd^{ti} \sum_{i \in I} \sum_{f \in F} Yd_{if}^{ti} Dd_{if} L_{if}^{ti} \tag{11}$$

The first set of constraints (13) are on the maximum capacity of the respective factory:

$$\sum_{w \in W} \sum_{tf \in TF} Ya_{fw}^{tf} \leq Pa_f Xa_f, \ \forall f \tag{12}$$

Constraints (14) on the maximum capacity of the respective warehouse:

$$\sum_{f \in F} \sum_{tf \in TF} Ya_{fw}^{tf} \leq Pb_w Xb_w, \ \forall w \tag{13}$$

Constraints (15) on the capacity of factory and warehouse:

$$\sum_{c \in C} \sum_{tw \in TW} Yb_{wc}^{tw} \leq \sum_{f \in F} \sum_{tf \in TF} Ya_{fw}^{tf}, \ \forall w \tag{14}$$

Constraints (16) and (16) on customer demand:

$$\sum_{w \in W} \sum_{tw \in TW} Yb_{wc}^{tw} \geq Q_c, \ \forall c \tag{15}$$

$$\sum_{i \in I} \sum_{tk \in TK} Yc_{ci}^{tk} \leq Q_c, \ \forall c \tag{16}$$

Constraint (18) on the maximum capacity of respective disassembly center:

$$\sum_{c \in C} \sum_{tk \in TK} Yc_{ci}^{tk} \leq Pd_i Xd_i, \ \forall i \tag{17}$$

Constraint (19) on the total number of demand of respective customer:

$$\sum_{i \in I} \sum_{tk \in TK} Yc_{ci}^{tk} \geq HdQ_c, \ \forall c \tag{18}$$

Constraint (20) on the product to be disposed of entered the respective disassembly center:

$$\sum_{f \in F} \sum_{ti \in TI} Yd_{if}^{ti} \geq Hr \sum_{c \in C} Yc_{ci}, \ \forall i \tag{19}$$

Constraint (21) on the maximum remanufacturing capacity of respective capacity:

$$\sum_{i \in I} \sum_{ti \in TI} Yd_{if}^{ti} \leq Pr_f Xa_f, \ \forall f \tag{20}$$

Non-negativity constraints (22) on decision variable:

$$Ya_{fw}^{tf}, Yb_{wc}^{tw}, Yc_{ci}^{tk}, Yd_{if}^{ti} \geq 0 \tag{21}$$

Binary number (23) explains the existence of facilities (factories, warehouses and disassembly center):

$$Xa_f, Xb_w, Xd_i \in \{0,1\} \tag{22}$$

## 4 Genetic Algorithm (GA)

The genetic algorithm relies on natural selection principles and evolutionary genetics. They combine the survival of the most effective string systems. The structure allows a randomized exchange of information to create a groundbreaking algorithm for optimality searching. Genetic Algorithms have been developed by Holland (1975), a population-based probabilistic and optimization technique. They are classified as global search heuristics and are also a class of evolutionary algorithms that use techniques triggered by evolutionary biologists such as inheritance, mutation, selection, and crossover (also known as recombination). They are implemented as a computer simulation of candidate solutions to the optimization problem, which communicates a better solution. Conventionally, solutions are represented in the 0s and 1s binary strings, but other encodings are not impossible. Evolution typically starts with a population of randomly generated individuals and happens in generations. In each generation, the fitness of each person in the population is measured, several individuals are selected from the current population (based on fitness) and changed (recombined and probably mutated) to create a new population. The new population is also used for the next iteration of the algorithm. The algorithm is normally completed when either a maximum number of generations has been created or an optimum fitness level has been reached for the population. If the algorithm has been terminated due to a maximum number of generations, an optimal solution may or may not have been found. For further information, please see a closed comprehensive example by Afshar-Nadjafi and Arani (2014).

## 5 Epsilon Constraint Method

The ε-constraint method was suggested by Haimes Y. V. et al. (1971). This method lessens the number of objectives of multi-objective problems. Only one of the objectives is considered as main objective functions and other objective functions are turned into constraints. So, the multi-objective problems are turned into single-objective problems, then it is called ε-constraint problem. By changing ε values with a proper interval, several Pareto optimal solutions can be found. The ideal points and nadir points can be determined and eliminated then by acquiring the full Pareto front. The set of ideal points and nadir points specifies the lower and upper bounds of Pareto optimal points.

## 6    Case Description

The planning horizon is 1 week long. The model includes two factories (plants), two warehouses, two customers (wholesalers), and two disassembly centers. The model is a closed-loop network. Fixed costs data, variable costs data, factory capacity data, transportation costs data, and demand of customers data have been collected from the Tsiakis and Shah (2016). Distance between stages data has been used from the Google maps (2020) and the rate of $CO_2$ released data has also been used from the Song et al. (2014).

## 7    Result Analysis

In **Table 2**, the decision variables of the mathematical model of the GSC have been obtained with the assistance of the Genetic Algorithm (GA) optimization technique in MATLAB R2017a software. The objective functions values are: $f_1 = £\ 1476934.685$ and $f_2 = 1885905.948(\text{kg})$.

**Table 4.** Genetic Algorithm's Output.

| Variables | Values | Variables | Values | Variables | Values |
|---|---|---|---|---|---|
| $Xa_1$ | 0.719 | $Xa_2$ | 0.673 | $Xb_1$ | 0.805 |
| $Xb_2$ | 1.304 | $Xd_1$ | 0.724 | $Xd_2$ | 0.758 |
| $Ya_{11}^1$ | 1629.999(units) | $Ya_{12}^1$ | 1629.999(units) | $Ya_{21}^2$ | 949.999(units) |
| $Ya_{22}^2$ | 949.999(units) | $Yb_{11}^1$ | 1629.999(units) | $Yb_{12}^1$ | 949.999(units) |
| $Yb_{21}^2$ | 1629.999(units) | $Yb_{22}^2$ | 949.999(units) | $Yc_{11}^1$ | 325.999(units) |
| $Yc_{12}^1$ | 325.999(units) | $Yc_{21}^2$ | 142.499(units) | $Yc_{22}^2$ | 142.499(units) |
| $Yd_{11}^1$ | 32.599(units) | $Yd_{12}^1$ | 14.249(units) | $Yd_{21}^1$ | 16.299(units) |
| | | $Yd_{22}^1$ | 7.125(units) | | |

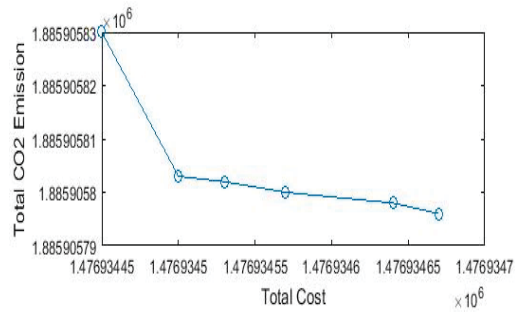

**Fig. 2.** Pareto Front of the Genetic Algorithm.

In this Pareto line, as the overall cost increases, the total CO2 emission decreases. But when the total cost decreases, the total emission of CO2 increases. There is a controversial relationship between them. There are some optimal points found in the Pareto front. Consequently, any company can use these optimal points for their purposes.

In **Table 3**, the decision variables of the mathematical model of the GSC have been obtained with the assistance of the ε-constraint technique in MATLAB R2017a software. The objective functions values are: $f_1 = £\ 1476877.1948$ and $f_2 = 1885570.948$(kg).

**Table 5.** ε-constraint Method Output.

| Variables | Values | Variables | Values | Variables | Values |
|---|---|---|---|---|---|
| $Xa_1$ | 0.718 | $Xa_2$ | 0.673 | $Xb_1$ | 0.805 |
| $Xb_2$ | 1.303 | $Xd_1$ | 0.725 | $Xd_2$ | 0.758 |
| $Ya^1_{11}$ | 1629.146(units) | $Ya^1_{12}$ | 1629.142(units) | $Ya^2_{21}$ | 949.655(units) |
| $Ya^2_{22}$ | 949.887(units) | $Yb^1_{11}$ | 1629.119(units) | $Yb^1_{12}$ | 950.718(units) |
| $Yb^2_{21}$ | 1629.111(units) | $Yb^2_{22}$ | 950(units) | $Yc^1_{11}$ | 326.105(units) |
| $Yc^1_{12}$ | 326.115(units) | $Yc^2_{21}$ | 142.732(units) | $Yc^2_{22}$ | 142.663(units) |
| $Yd^1_{11}$ | 32.665(units) | $Yd^1_{12}$ | 14.333(units) | $Yd^2_{21}$ | 16.377(units) |
| | | $Yd^2_{22}$ | 7.174(units) | | |

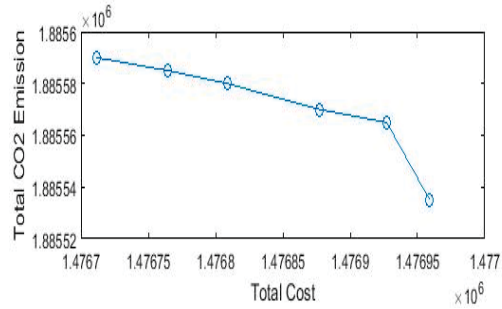

**Fig. 3.** Pareto Front of the ε-Constraint Method.

In this Pareto front, only the dominating points are used. When the total cost increases, the total $CO_2$ emission decreases. Otherwise, when the total cost decreases, therefore, the total $CO_2$ emission increases. There are some optimal points found in this Pareto front by solving this method.

## 8    Conclusion and Future Research

This study is attached to the problem context. There is a controversial relationship between the two objectives of total cost and total $CO_2$ emission. For examining the optimal decision of the objectives, a mathematical model is formulated. The genetic algorithm has become an effective technique in solving many critical problems and has applied to practical problems. Epsilon Constraint Method is also a popular method for solving many practical problems. The mathematical model is solved by the Genetic algorithm and Epsilon Constraint Method, individually. The comparison of these two methods is one of the contributions of this study. Secondary data is used to examine to verify the mathematical formulation. For transport choices, road transport is selected as a single mode of transport for all stages of the supply chain. The cost of transport depends on the variance in the unit cost of transport and on the total amount of product units flowing. The cost of transport, therefore, has a major impact on the overall cost. The overall emission of $CO_2$ also depends heavily on the total distance. The distance (km) × kg is proportional to the overall $CO_2$ emission. Therefore, the overall emission varies with distance variability. Other factors have not had a significant impact on total $CO_2$ emissions. Thus, the total cost and total $CO_2$ emission decrease when unit transport costs and distances are reduced. The genetic algorithm and ε-constraint method have been compared. But closely the same results have been obtained by these two methods.

Currently, the GSC area increases scholars' and governments' awareness because of the recent environmental issue. A company thinks mainly the financial fact and concerns about environmental issues. So, this paper suggests a mathematical model of the GSC, and its methodology solve with a multi-objective approach. Total cost and total $CO_2$ emission are two main criteria to evaluate which are selected to construct a mixed-integer linear programming model. A Pareto front represents the two objectives of total cost and total emission of $CO_2$. One closed-loop network design is accompanied by the problem statement. The assumptions lead to a more realistic model. The proposed genetic algorithm and the ε-constraint method are intended to minimize total costs along with total $CO_2$ emissions. Over the last few years, researchers have proposed several models and studies for the design of the green supply chain with varying degrees of complexity. This research is an attempt to bridge between the total cost of supply and total $CO_2$ emissions and to find the optimum point for these two objectives.

The mathematical model of the GSC design is confined to its services. So, the model will be extended for increasing its services. Then the model is more applicable in complex situations. Further analysis can then be tackled by changing the model scenario. Other optimization techniques such as fuzzy multi-objective linear programming (FMLOP) (Momeni Tabar et al. 2014), Discrete Event Simulation (DES) (Arani et al. 2020) can be used to solve the formulated model. Then compare the optimization technique which leads us to draw a comprehensive comparison. For the delivery of products, different modes of transport can be used. Additionally, compare the different modes of transport and find the right one for this model. These recommendations could shed light on future research.